\documentclass[12pt]{article}

\usepackage{amsmath, amsfonts, amscd}
\usepackage[all]{xy}

\newtheorem{thm}{Theorem}[section]
  \newtheorem{prop}[thm]{Proposition}

      \newtheorem{cor}[thm]{Corollary}
         \newtheorem{lem}[thm]{Lemma}
   
\newtheorem{ex}[thm]{Example}

\newtheorem{rem}[thm]{Remark}

\newcommand{\bZ}{\mathbb{Z}}
\newcommand{\bR}{\mathbb{R}}

\newcommand{\ul}{\underline}
\newcommand{\cQ}{\mathcal{Q}}
\newcommand{\La}{\Lambda}
\begin{document}

\begin{center}{\bf\Large Ext-finite modules
for weakly symmetric algebras with radical cube zero}

\medskip

Karin Erdmann

\vspace*{1cm}

Dedicated to the memory of Laci Kov\'{a}cs.

\begin{abstract}
\noindent
Assume $A$ is weakly symmetric, indecomposable, with
radical cube zero and radical square non-zero. We show that such algebra 
of wild representation type does not have
a non-projective module $M$ whose ext algebra is finite-dimensional.
This gives a complete classification  weakly symmetric indecomposable
algebras which have  a non-projective module whose
ext algebra is finite-dimensional.

\medskip
\noindent
\textbf{MR Subject classification:} 16E40, 16G10, (16E05, 15A24, 33C45)

Chebyshev polynomials. 

\medskip

\end{abstract}
\end{center}

\parindent=0pt
\baselineskip=16pt
\textwidth=15cm
\textheight=22cm

\section{Introduction}

Assume $A$ is a finite-dimensional algebra over a field $K$. We say that an
$A$-module $M$ is ext finite if there is some $n\geq 0$ such that 
${\rm Ext}^k_A(M, M)=0$ for $k>n$. 

If $A=KG$, the group algebra of a finite group, then any ext finite module
is projective (this may be found in Chapter 5 of \cite{Be}). On the other hand, there is a four-dimensional selfinjective algebra
which has non-projective ext finite modules, first described in 
\cite{Sch}. This algebra is known as q-exterior algebra, see Section 4. 
If a selfinjective algebra $A$ has a non-projective ext finite module 
there is no support variety theory for $A$-modules via Hochschild cohomology.
This follows from Corollary 2.3  in \cite{EHSST}, it shows that the finite generation
conditions (Fg1, 2) in \cite{EHSST} (and (Fg) of \cite{So}) must fail. 
That is, existence of ext-finite non-projective modules
gives information about action of the Hochschild cohomology on 
ext algebras of modules. 

\newpage

There is also the 'generalized Auslander-Reiten 
condition', GARC, which has been introduced in \cite{AR} in the context of homological conjectures, which has attracted a lot of interest,  see for example \cite{CH, CT, Jo,   DP1, DP2}. 
The condition  GARC for a ring $R$ states:
 
{\it  If $M$ is an $R$-module and there is some $n\geq 0$ such that
${\rm Ext}^k_R(M, M\oplus R)=0$ for $k>n$, then $M$ has projective dimension
at most $n$.}

The four-dimensional local algebra mentioned above does not satisfy GARC, 
there are even  counterexamples
with $n=1$, see Section 4.
It is not known whether there is a ring $R$ which has a counterexample with $n=0$. 

If $R= A$ and $A$ is a selfinjective finite-dimensional algebra
then GARC states that  
\ {\it  any ext-finite module is projective}.

\medskip

The four-dimension algebras which have non-projective ext finite modules
belong to the class of weakly symmetric algebras with  radical cube zero. 
These algebras have been studied in \cite{Be2, ESo}. In particular it is 
understood when such algebra does not satisfy the (Fg) condition:

\ \ 
Assume $A$ is weakly symmetric with radical cube zero (and radical square non-zero). Assume also
$A$ is indecomposable.
Let $E$ be the matrix with entries $\dim {\rm Ext}^1(S_i, S_j)$ where 
$S_1, S_2, \ldots, S_r$ are the simple $A$-modules. Then $E$ is a symmetric
matrix, so it has real eigenvalues. The  largest eigenvalue 
$\lambda$ say, occurs with multiplicity one, and has a positive eigenvector,
this is the Perron-Frobenius Theorem. 
It is proved in \cite{ESo} that $A$ does not satisfy (Fg) if and only if
either $\lambda > 2$, or else $A$ is Morita equivalent to a four-dimensional
local algebra as above, or else to a 'Double Nakayama algebra' (see section 4), where  in both 
cases there is  a deformation parameter which is not a root of unity.

\bigskip

These Double Nakayama algebras  also  have ext-finite
non-projective modules; this is probably known, 
we will give a proof in Section 4. 

Our main result   shows that a
weakly symmetric algebra with radical cube zero and $\lambda>2$ does not have 
ext-finite non-projective modules. With this, we get the following.

\bigskip
\begin{thm} Assume $A$ is a weakly symmetric indecomposable algebra
over an algebraically closed field, with 
$J^3=0\neq J^2$. Then $A$ has an ext-finite non-projective module if and only if
$\lambda=2$ and $A$ is Morita equivalent to either a four dimensional q-exterior algebra, 
or a Double Nakayama algebra, where in both cases the deformation parameter
is not a root of unity.
\end{thm}

\medskip

It follows that
existence of  ext finite non-projective modules
is not equivalent with failure of (Fg).

\medskip

The Theorem remains true for an arbitrary field if one takes for 
$A$ an algebra defined by quiver and relations.

Section 2 contains the relevant background. In Section 3 
we prove  the main new part of the Theorem,  and in Section 4 we describe ext-finite non-projective modules for 
the algebras for which $\lambda=2$.

We work with finite-dimensional left $A$-modules, and if $M, N$ are 
such $A$-modules then we write ${\rm Hom}(M, N)$ instead of ${\rm Hom}_A(M, N)$
and similarly ${\rm Ext}^k(M, N)$ means ${\rm Ext}^k_A(M, N)$.

\vspace*{2cm}

\section{Preliminaries}

{\bf 2.1 } 
We assume throughout that $A$ is a finite-dimensional weakly symmetric 
algebra over an algebraically closed  field $K$, and we assume $A$ is indecomposable. This is no restriction since we will focus on 
indecomposable modules. Suppose
$M$ is a finite-dimensional $A$-module. Then 
${\rm rad}(M)$ is the submodule of $M$ such that $M/{\rm rad}(M)$ is the
largest semisimple factor module of $M$, sometimes called 'top' of $M$. 
The module ${\rm rad}(M)$ is equal to $JM$ where $J$ is the radical of $A$. 
The socle of $M$, denoted by ${\rm soc}(M)$, is the largest semisimple submodule
of $M$. 

\medskip

{\bf 2.2 } A finite-dimensional $A$-module $M$ has a projective cover, that is
there is a surjective map $\pi_M: P\to M$ where 
$P$ is projective, and $P/{\rm rad}(P)\cong M/{\rm rad}(M)$. The kernel of $\pi_M$ is unique up to isomorphism, and is denoted by $\Omega(M)$. 
Repeatedly taking projective covers gives a minimal projective resolution of $M$, 
$$\ldots \to P_m \stackrel{d_m}\to P_{m-1} \to \ldots \stackrel{d_1}\to P_0 
\stackrel{d_0}\to M\to 0
$$
where $d_0= \pi_M$ and $d_m$ is a projective cover of $\Omega^m(M)$ for $m\geq 1$. 
If $A$ is selfinjective and $M$ is indecomposable and non-porojective
then also $\Omega(M)$ is indecomposable and non-projective. In fact, $\Omega$ induces an equivalence
of the stable module category of $A$.

{\bf 2.3 } \ 
We assume $A$ is weakly symmetric. This means that $A$ is selfinjective, and 
any indecomposable projective
module has a simple socle, with ${\rm soc}(P) \cong P/{\rm rad}(P)$. Hence
for any  simple module, its projective cover is also its injective hull. 
This implies also that for any non-projective 
indecomposable $A$-module $M$ we have
that $M/{\rm rad}(M$) is isomorphic to ${\rm soc}\ \Omega(M)$.

Let $S_1, S_2, \ldots, S_r$ be the simple $A$-modules, and let
$P_i$ be the projective cover of $S_i$ for $1\leq i\leq r$. We assume $J^3=0$ but $J^2\neq 0$. 
If so, then every indecomposable projective $P_i$ must have
radical length three; this is well known (and it is easy to see, recalling
that we assume $A$ to be indecomposable). 
So we have $P_i/{\rm rad}(P_i) \cong S_i \cong {\rm soc}(P_i)$ and
${\rm rad}(P_i)/{\rm soc}(P_i)$ is semisimple and non-zero. So we can write
$${\rm rad}(P_i)/{\rm soc}(P_i) \cong \oplus _{j=1}^r S_j^{d_{ij}}
$$ 
where $d_{ij}\geq 0$ and not all $d_{ij}$ are zero. It is also true that
for all $i, j$ we have $d_{ij} = d_{ji}$. We will give the proof in 2.6
below.

 This is  sufficient information to compute dimensions
of $\Omega$-translates of $M$. The crucial property is the following, which is well-known. 
For convenience we give the proof.

\medskip

\begin{lem} \label{lem:socrad} \
Assume $M$ is a module with no simple or projective summands. Then
${\rm soc}(M) = {\rm rad}(M)$.
\end{lem}

{\it Proof } Since $M$ has no projective (hence injective ) summand
it has radical length $\leq 2$. Therefore ${\rm rad}(M)$ is annihilated
by $J$ and hence is contained in ${\rm soc}(M)$. The socle of $M$ is
semisimple, hence 
${\rm soc}(M) = {\rm rad}(M) \oplus C$ for some submodule $C$ of ${\rm soc} (M)$. We must show that $C=0$.

Let $\pi: M\to M/{\rm rad}(M)$ be the canonical surjection, then $\pi(C)$ is
isomorphic to $C$, write $C' = \pi(C)$.

The module $M/{\rm rad}(M)$ is simisimple, so we can write
$M/{\rm rad}(M) = C'\oplus G$ for some semisimple module $G$. Let
$\tilde{G}$ be the submodule of $M$ containing ${\rm rad}(M)$ such
that $\tilde{G}/{\rm rad}(M) = G$.

Then we have that $\tilde{G}\cap C=0$ and $M= \tilde{G} + C$:
Namely if $x\in \tilde{G}\cap C$ then $x+{\rm rad}(M) \in G\cap C'=0$ and therefore $x\in {\rm rad}(M)$, and then it is in the intersection
of ${\rm rad}(M)$ with $C$ and is zero.
Furthermore, we  have $M/{\rm rad}(M) = G+C'$which imples that $M= \tilde{G} + C$. So if $C\neq 0$ then it is a semisimple summand of $M$, and by the assumption $C=0$.
$\Box$

\bigskip

{\bf 2.4 }\label{dimvec} \ Let $M$ be a module such that ${\rm soc}(M) = {\rm rad}(M)$, both
the socle of $M$ and $M/{\rm rad}(M)$ are semisimple.
We write $\ul{\dim}({\rm soc}(M) = \ul{s} = (s_1, s_2, \ldots, s_n)$ where
$s_i$ is the multiplicity of $S_i$ in ${\rm soc}(M)$, and similarly
we write\\ 
$\ul{\dim}(M/{\rm rad} M)= \ul{t} = (t_1, t_2, \ldots, t_n)$
where $t_i$ is the multiplicity of $S_i$ in 
$M/{\rm rad} \ M$. 

Then we define the  'dimension vector' for $M$ to be 
$$\ul{\dim} (M) := (\ul{t} \ |\  \ul{s})$$
The usual dimension vector would be $\ul{t} + \ul{s}$.

\begin{lem}\label{lem:dim} 
Let $X$ be the $2n\times 2n$ matrix which in block form is
given by 
$$X = \left(\begin{matrix} E & -I_n \cr I_n & 0 \end{matrix}\right).
$$
Assume  $M$ has no simple  or projective summands, and 
$\Omega(M)$ is not simple. Then 
$$\ul{\dim} \Omega(M)^T  = X\ul{\dim}(M)^T.
$$
\end{lem}

{\it Proof } Consider the projective cover of $M$,
$$0 \to \Omega(M) \to P_M \to M\to 0
$$
Then $P_M \cong \oplus_{i=1}^n t_iP_i$ since $P_M/{\rm rad}(P_M)$ must be
isomorphic to $M/{\rm rad}(M)$. 
Since $M$ has no projective ( hence injective) summands, the socle
of $\Omega(M)$ is isomorphic to ${\rm soc}(P_M)$ that is $\oplus t_iS_i$. 

As well, since $\Omega(M)$ has no simple or projective summand, we know
${\rm soc}\ \Omega(M) = {\rm rad}\ \Omega(M)$, 

\bigskip

Factoring out the socle of $\Omega(M)$ we get a short exact sequence
$$0 \to \Omega(M)/{\rm soc} \Omega(M) \to P_M/{\rm soc}(P_M) 
\to M\to 0
$$
If we restrict this to the radical of
$P_M/{\rm soc}(P_M)$ then we get a split exact sequence,
$$0 \to \Omega(M)/{\rm soc} \Omega(M)  \to {\rm rad}(P_M)/{\rm soc}(P_M)
\to {\rm soc}(M)\to 0
$$
Hence the dimension vector of  $\Omega(M)/{\rm soc} \Omega(M)$ 
is equal to 
$$E\ul{t}^T - \ul{s}^T$$
as required.
$\Box$

\medskip
This would still be true if $\Omega(M)$ is simple. 
Since we want to iterate the calculation, we exclude this.

\bigskip

{\bf 2.5 } If none of the the modules
$\Omega^r(M)$ for $r=1, 2, \ldots, k+1$ is simple, it follows  that
the dimension vector of $\Omega^k(M)$ is equal to
$X^k\ul{\dim}(M)^T$. 
The matrix 
$X^k$ is of the form
$$X^{k} = \left(\begin{matrix} f_{k}(E) & -f_{k-1}(E)\cr f_{k-1}(E) & -f_{k-2}(E)
\end{matrix}\right).
$$
Here $f_m(x)$ is the $m$-th Chebyshev polynomial, given by
$$ f_0(x)=1, \ \  f_1(x) = x, \ \  
f_{k}(x) = xf_{k-1}(x) - f_{k-2}(x)  \ (k\geq 2).
$$
The polynomial $f_k(x)$ is the characteristic polynomial of the
$k\times k$ incidence matrix of the Dynkin diagram of type $A_k$, that is it has
entries
 $a_{i,i\pm 1}= 1$ and $a_{ij}=0$ otherwise.

As well, $f_k(x) = U_k(x/2)$ where $U_k(x)$ is a version of a 
Chebyshev polynomial of the second kind.
These polynomials are studied extensively in numerical mathematics, see for example
\cite{Ri}.

\bigskip

{\bf 2.6 } We recall that if $S$ is a simple module
then for any $k\geq 1$ and for any module $M$ we have
$${\rm Ext}^k(M, S) = {\rm Hom}(\Omega^k(M), S)$$

We give the argument.
Take the exact sequence 
$$0\to \Omega^k(M)\to P_{k-1}\to \Omega^{k-1}(M)\to 0
$$
and apply ${\rm Hom}(-, S)$.
If $\pi: P_{k-1}\to S$ then clearly this restricts to the zero map $\Omega^{k}(M)\to S$.  Hence
the inclusion ${\rm Hom}(\Omega^{k-1}(M), S)\to {\rm Hom}(P_{k-1}, S)$ is an isomorphism.
Therefore ${\rm Hom} (\Omega^k(M), S)\cong {\rm Ext}^1(\Omega^{k-1}(M), S)$
which is isomorphic to ${\rm Ext}^k(M, S)$.

We claim that  $d_{ij} = d_{ji}$, which shows that the matrix $E$ is symmetric: 
We may assume $i\neq j$. Then since $P_j$ is the projective cover of $S_j$ we have
$$d_{ij} = \dim {\rm Hom}(P_j, P_i)
$$
But $P_i$ is also the injective hull of $S_i$, so the dimension is also equal to the number
of times $S_i$ occurs in $P_j$, which is equal to $d_{ji}$.

\vspace*{1cm}

\section{The main result}

Assume that $A$ is weakly symmetric with $J^3=0 \neq J^2$ and 
$A$ has an ext finite non-projective module, then there is such $M$ which
is indecomposable. We will analyse the
dimension vectors of the modules $\Omega^r(M)$ for large $r$.

We may assume that $\Omega^r(M)$ is not simple for $r\geq 0$: Namely at most
finitely many of the $\Omega^k(M)$ can be simple, since otherwise it would
follow that $M$ is periodic, but then $M$ would not be ext finite. 
So there is some $m$ such that for $k\geq m$ none of the modules  $\Omega^k(M)$ is simple.
We replace $M$ by $\Omega^m(M)$.   With $M$, also 
$\Omega^m(M)$ is ext-finite and not projective, recall that 
 $\Omega$ induces an equivalence of
the stable category.

\bigskip

\begin{lem} Assume ${\rm Ext}^k(M, M)=0$ for $k>n$. 
Let $(\ul{t}\ |\ \ul{s})$ be the dimension vector of $M$ and 
 $(\ul{t}^{(k+1)}\ |\  \ul{s}^{(k+1)})$ be
the dimension vector of $\Omega^{k+1}(M)$. Then for all $k>n$ we have
$$(\ul{s}\ |\ -\ul{t})\cdot (\ul{t}^{(k+1)}\ |\ \ul{s}^{(k+1)}) = 0.
$$
\end{lem}

{\it Proof } By the assumption, and by \ref{lem:socrad},   ${\rm soc}(M) = 
JM$, so  we have 
a short exact sequence
$$0\to M_2= {\rm soc}(M) \to M \to M_1 = M/JM \to 0
\leqno{(*)}$$
where $M_1$ and $M_2$ are semisimple.
Write $M_1 = \oplus_i t_iS_i$ and $M_2 = \oplus_i s_iS_i$. Then 
${\rm dim }(M) = (\ul{t}\ |\ \ul{s})$.

We apply the functor ${\rm Hom}(M, -)$ to (*) which gives
the long
exact sequence of homology. Part of this is 
$$\ldots \to {\rm Ext}^k(M,M) \to {\rm Ext}^k(M, M_1)
\to {\rm Ext}^{k+1}(M, M_2) \to {\rm Ext}^{k+1}(M, M) \to \ldots .
$$
Consider 
${\rm Ext}^k(M, M_1)$,  this is
isomorphic to $\oplus_i t_i {\rm Ext}^k(M, S_i)$, and
${\rm Ext}^k(M, S_i) \cong {\rm Hom}(\Omega^k(M), S_i)$ (see 2.6).
This has dimension
$$\sum_i t_it_i^{(k)}.
$$
Similarly ${\rm Ext}^{k+1}(M, M_2)$ has 
dimension
$$\sum_i s_it_i^{(k+1)}.
$$
By exactness we get for $k>n$ that
${\rm Ext}^k(M, M_1) \cong {\rm Ext}^{k+1}(M, M_2).
$
Equating  dimensions,  we have
$$\sum_i t_it_i^{(k)} = \sum_i s_it_i^{(k+1)}.
$$
By 2.3 we know that $\ul{t}^{(k)} = \ul{s}^{(k+1)}$. Using this, and
rewriting the last identity we get the claim.
$\Box$

\vspace*{1cm}

{\bf 3.2 } We analyse $(\ul{s}\ |\ -\ul{t})\cdot (\ul{t}^{(k+1)}\ |\ \ul{s}^{(k+1)})$, 
which is equal to
$$(\underline{s}\ | \ -\underline{t})X^{k+1}(\underline{t}\ | \ \underline{s})^T
\leqno{(1k)}
$$
for $k > n$. We substitute $X^{k+1}$ and expand, then (1k) becomes
$$\underline{s}f_{k+1}(E)\underline{t}^T - \underline{t}f_k(E)\underline{t}^T
- \underline{s}f_k(E)\underline{s}^T + \underline{t}f_{k-1}(E)\underline{s}^T
\leqno{(2k)}.
$$
The matrix $f_{k-1}(E)$ is symmetric, so we can intechange $\underline{t}$ and
$\underline{s}$ in the last term. Then using the recurrence relation for
the Chebyshev polynomials,
$$f_{k+1}(x) = xf_k(x) - f_{k-1}(x),$$
the expression (2k) becomes
$$\underline{s}Ef_k(E)\underline{t}^T - \underline{t}f_k(E)\underline{t}^T 
- \underline{s}f_k(E)\underline{s}^T.
\leqno{(3k)}
$$
Since $E$ is real symmetric, there is an orthogonal matrix $R$
such that $R^TER=D$, a diagonal matrix. We substitute $E=RDR^T$, and we
set $\underline{\alpha}:= \underline{s}R$ and
$\underline{\beta}:= \underline{t}R$.
With this, noting also that $Rf_k(E)R^T = f_k(RER^T)$, expression (3k) becomes
$$\underline{\alpha} Df_{k}(D)\underline{\beta}^T - \underline{\beta}f_k(D)
\underline{\beta}^T - \underline{\alpha}f_k(D)\underline{\alpha}^T. 
\leqno{(4k)}$$
The matrices involved are diagonal, let $\lambda_1, \ldots, \lambda_r$ be
the eigenvalues of $D$. Then (4k) is equal to
$$\sum_{i=1}^r(\alpha_i\beta_i \lambda_i - \beta_i^2 - \alpha_i^2)f_k(\lambda_i).$$
If we denote the distinct eigenvalues of $D$ 
by $\mu_1, \ldots, \mu_m$ then this
becomes
$$\sum_{j=1}^m (\sum_{\lambda_i=\mu_j} 
\alpha_i\beta_i \lambda_i - \beta_i^2 - \alpha_i^2)f_k(\mu_j).
\leqno{(5k)}$$

Then Lemma 3.1 shows that (5k) is zero for all $k>n$.
The coefficients $c_j:= (\sum_{\lambda_i=\mu_j} 
\alpha_i\beta_i \lambda_i - \beta_i^2 - \alpha_i^2)$
do not depend on $k$. We take any $m$ of these equations
for $k>n$ and write them in matrix form. That is, consider a matrix 
$$C:= \left(\begin{matrix} f_N(\mu_1)& f_N(\mu_2) & \ldots & f_N(\mu_m)\cr
f_{N+i_1}(\mu_1) & f_{N+i_1}(\mu_2) & \ldots & f_{N+i_1}(\mu_m)\cr
\ldots &&&\cr
f_{N+i_{m-1}}(\mu_1) & f_{N+i_{m-1}}(\mu_2) & \ldots & f_{N+i_{m-1}}(\mu_m)
\end{matrix}\right)
$$
with $N>n$ and $0<i_1<i_2 < \ldots < i_{m-1}$. Then for any such $C$ we have
$$C\left(\begin{matrix}c_1\cr c_2\cr \vdots \cr c_m\end{matrix}\right)=0.
$$

\begin{lem} There are $N>n$ and $0<i_1<i_2< \ldots < i_m$ such that
$C$ is non-singular.
\end{lem}

{\it Proof } 
\ We use induction on $m$. Assume first $m=1$. 
We have  $f_0(\lambda_1)=1\neq 0$.
Whenever $f_{u-1}(\lambda_1)\neq 0$ and $f_u(\lambda_1)=0$ then 
$f_{u+1}(\lambda_1) = -f_{u-1}(\lambda_1)\neq 0$. 
So at worst every second of the values can be zero.

Now we fix some $N>n$ such that $f_N(\lambda_1)\neq 0$. 
Consider the matrix with rows
$$R_{N+j}:= 
(f_{N+j}(\mu_1) ,  f_{N+j}(\mu_2) , \ldots , f_{N+j}(\mu_m))
$$
for $j=0, 1, 2, \ldots, k$ and $k$ large, $k>m+2$.
We replace $R_{N+k}$ by $R_{N+k} + R_{N+k-2} - \mu_1R_{N+k-1}$ and obtain as the new last row
$$[0, (\mu_2-\mu_1)f_{N+k-2}(\mu_2), \ldots, (\mu_m-\mu_1)f_{N+k-2}(\mu_m)]
$$
Similarly we replace $R_{N+k-1}$, and so on. This process ends when
row $R_{N+2}$ has become
$$[0, (\mu_2-\mu_1)f_{N+1}(\mu_2), \ldots, (\mu_m-\mu_1)f_{N+1}(\mu_m)]
$$
By construction, $f_N(\mu_1)\neq 0$,  and
we take the row of $f_N(\mu_i)$  as the first row of our required submatrix.

We apply the inductive hypothesis to the matrix consisting of $R_{N+2}, \ldots, R_{N+k}$
 omitting
the first column. 
Note that from each column we can take a non-zero scalar factor
$\mu_i-\mu_1$. 
The remaining matrix has the same shape again with $m-1$ columns. 
So by the inductive hypothesis it has $m-1$ rows which
form a non-singular submatrix.
$\Box$

\bigskip

\begin{ex}
\normalfont The roots of $f_r(x)$ are precisely the
eigenvalues of the $r\times r$ matrix $E$ with $e_{i, i\pm 1}=1$ and $e_{ij}=0$ otherwise (see 2.5). By the Cayley-Hamilton theorem we know  that  $f_r(E)=0$.
In \cite{ESch}  it is proved that the sequence of matrices $(f_m(E))$ is periodic.
In fact one can see from the proof there that there are $r$ successive
rows which are linearly independent, but there are rows of zeros.

For example $r=2$, then the eigenvalues are $\pm 1$ and the rows are
$$\begin{matrix} 1 & -1\cr 0 & 0 \cr -1 & 1\cr -1 & -1 \cr 0 & 0 \cr 1 & 1 \cr 1 & -1 \cr 0 & 0 \cr \ldots & \end{matrix}.
$$
\end{ex}

\begin{cor} If (1k) is zero for all $k>n$ then for all $j$ with 
$1\leq j\leq m$ we have
$$\sum_{\lambda_i=\mu_j} (\alpha_i\beta_i\lambda_i - \beta_i^2 - \alpha_i^2)=0.
$$
\end{cor}

This follows from the previous Lemma.

\medskip

Let $\lambda_1 = \lambda$, the largest eigenvalue of $E$. 
We assume $A$ is indecomposable and therefore $E$ is an irreducible matrix.
Therefore $\lambda$ has multiplicity one as an eigenvalue of $E$, and there is a real
eigenvector $\ul{v}$ with $v_i>0$ for all $i$. We may take it as a
unit vector, and then $\ul{v}^T$ 
is the first column of $R$ where $R^TER =D$.  
Recall $\ul{\alpha} = \ul{s}R$ and $\ul{\beta} = \ul{t}R$. These have
first components
$$\alpha_1 = \sum_i s_iv_i, \ \ \beta_1 = \sum_i t_iv_i.
$$
Since $\ul{s}$ and $\ul{t}$ are non-zero in $\bZ^r_{\geq 0}$ it follows
that $\alpha_1$ and $\beta_1$ are positive. 
Because $\lambda$ has multiplicity one, the sum in 3.4 for $\lambda$
has only one term, and we deduce:

 \begin{cor} The numbers $\alpha_1/\beta_1$ and $\beta_1/\alpha_1$ are roots of the equation
 $$X^2 - \lambda X +1$$
 \end{cor}

 We may  do the same calculation with $\Omega^m(M)$
instead of $M$ for any $m\geq -1$, denote the  corresponding  numbers by 
$\beta^{(m)}_1$ and $\alpha^{(m)}_1$. For any such $m$, the two
quotients must therefore be 
 roots of the above  quadratic equation, that is
$$\frac{\beta_1^{(m)}}{\alpha_1^{(m)}} + \frac{\alpha_1^{(m)}}{\beta_1^{(m)}}
= \lambda \leqno{(**)}
$$ 
We can say more.

 \begin{lem} 
We have $\beta^{(m+1)}_1 = \lambda \beta_1^{(m)} - \alpha_1^{(m)}$. 
 \end{lem}

 \bigskip
 
 {\it Proof }  For the proof, it suffices to take $m=0$. We have $\ul{t}^{(1)} = E\ul{t}^T - \ul{s}$ and therefore 
 $$\ul{t}^{(1)}_i = (E\ul{t}^T)_i- s_i$$
 Now $(E\ul{t}^T)_i = \sum_{k=1}^r e_{ik}t_k = \sum_{k=1}^r e_{ki}t_k$ (recall $E$ is symmetric). Then 
 $$\beta^{(1)}_1 + \alpha^{(0)}_1 = \sum_{i=1}^n (E\ul{t}^T)_iv_i.
$$ 
We substitute and change the order of summation and get that this equal to
$$\sum_{k=1}^r(\sum_{i=1}^r e_{ki}v_i) t_k.
$$
The 
coefficient of $t_k$ is  the k-th entry of $E\ul{v}^T = \lambda \ul{v}$, which is $\lambda v_k$. So we get
$$\beta_1^{(1)} + \alpha_1^{(1)} = \lambda\sum_k v_kt_k = \lambda\beta^{(0)}_1,
$$
as stated.
$\Box$

\bigskip

 \begin{prop} If $M$ is ext-finite then $\alpha_1 = \beta_1$. 
In particular $\lambda = 2$. \end{prop}

\bigskip

{\it Proof }  (1) First we claim that
$\alpha_1^{(m)}/\beta_1^{(m)} = \alpha_1^{(m+1)}/\beta_1^{(m+1)}$.

By  Lemma 3.6, and since 
$\alpha_1^{(m+1)} = \beta_1^{(m)}$ (recall $\ul{s}^{(m+1)} = \ul{t}^{(m)}$), 
 we have 
$$\frac{\beta_1^{(m+1)}}{\alpha_1^{(m+1)}} =  \lambda - \frac{\alpha_1^{(m)} }{\beta_1^{(m)}}
$$
Using also (**)  we deduce 
$$\frac{\beta_1^{(m+1)}}{\alpha_1^{(m+1)}} + \frac{\alpha_1^{(m+1)}}{\beta_1^{(m+1)}} = 
\lambda = \lambda + \frac{\alpha_1^{(m+1)}}{\beta_1^{(m+1)}} - \frac{\alpha_1^{(m)}}{\beta_1^{(m)}}
$$
and hence  the claim follows.

\medskip
The set of positive numbers $\{ \alpha_1^{(m)}, m\geq -1 \}$ is bounded below,
 and it is a discrete subset of $\bR$, therefore it has a minimum. 
That is, 
we may  choose $M$ in its $\Omega$-orbit so that the number $\alpha_1^{(1)} \leq \alpha_1^{(m)}$ for all $m \geq -1$.

Then $\beta_1^{(1)} = \alpha_1^{(2)} \geq \alpha_1^{(1)} = \beta_1^{(0)}$ and $\alpha_1^{(1)}\leq \alpha_1^{(0)}$. It follows that

$$\frac{\alpha_1^{(0)}}{\beta_1^{(0)}} = \frac{\alpha_1^{(0)}}{\alpha_1^{(1)}} \geq 1
$$
$$\frac{\alpha_1^{(1)}}{\beta_1^{(1)}} = \frac{\alpha_1^{(1)}}{\alpha_1^{(2)}} \leq 1 
$$
and hence the fractions must be equal to 1.

So the quadratic equation has one root equal to 1. The  product of the
roots is 1, so both roots are equal to 1 and then $\lambda =2$.
$\Box$

\bigskip

We have proved  that for $\lambda \neq 2$, the algebra has no ext-finite modules.

\bigskip

\begin{rem}\normalfont
Assume $\lambda =2$. For the algebras without (Fg), (which are of type 
$\tilde{A}$, or local,)
the vector $\ul{v}$
is  a multiple of 
$(1, 1, \ldots, 1)$, and if 
 $\alpha_1^{(m)}=\beta_1^{(m)}$ for all $m$ then the 
socle and the top of any $\Omega$-translate of $M$ have the 
same dimension. 
So $M$ has even dimension, and it follows that $M$ cannot be an $\Omega$ translate of a simple module. Namely the  $\Omega$ translates of simple modules have odd dimensions for these algebras.
\end{rem}

\bigskip

\vspace*{1cm}

\section{Algebras where $\lambda=2$}

If $A$ has an ext-finite non-projective module then (Fg) does not hold. By \cite{ESo} when $\lambda=2$, the algebra is
Morita equivalent to either the q-exterior algebra, or to an algebra of type $\widetilde{A}_n$, which we call a 
Double Nakayama algebra. In both cases, there is a deformation parameter which is not a root of $1$ (and non-zero).

\subsection{The q-exterior algebra} 

Let $\La = \La(q) =
K\langle x, y\rangle / (x^2, y^2, xy+qyx)$ and $0\neq q\in K$.
We assume that $q$ is not a root of unity. 

It was discovered by R. Schulz, already some years ago, that
this algebra has ext finite non-projective modules, see section 4 in \cite{Sch}.

For $0\neq \lambda \in K$ we define a $\La$-module $M= C(\lambda)$ 
as follows. It is 2-dimensional and $x, y$ act by
$$x\mapsto \left(\begin{matrix} 0 & 1\cr 0 & 0\end{matrix}\right),  \ \ y\mapsto \left(\begin{matrix} 0 & \lambda\cr 0 & 0\end{matrix}\right), 
$$
It is clearly indecomposable and not projective, and it is easy to check that
 $C(\lambda)\cong C(\mu)$ only if $\lambda=\mu$. 
We  construct $C(\lambda)$ as the submodule of $\La$ generated by 
 $\zeta = -\lambda qx + y\in \La$, and take basis $\zeta, x\zeta$.
 \medskip

 \begin{lem} We have $\Omega^m(C(\lambda)) \cong C(q^{-m}\lambda)$ for $m \geq 0$.
 \end{lem}

 {\it Proof }  
We find $\Omega(M) = \{ z\in \La: z\zeta = 0\} = \La(y-\lambda x)$
and if $\zeta_1 = y-\lambda x$ then $y\zeta_1 = \lambda q^{-1}x\zeta_1$.  That is,
$\Omega(M) \cong  C(\lambda q^{-1})$, and the statement follows.

\bigskip

For convenience we give a proof showing that the module 
$C(\lambda)$ is ext-finite.

\begin{lem} If $\mu\in K$ and $\mu\neq \lambda$ or $q\lambda$ then ${\rm Ext}^1_(C(\mu), C(\lambda))=0$. 
\end{lem}

\bigskip

{\it Proof }  A projective cover of $C(\mu)$ is  of the form
$$0\to C(\mu q^{-1}) \to \La \to C(\mu)\to 0
$$
Applying ${\rm Hom}(-, C(\lambda))$ gives a four term exact sequence
$$0 \to {\rm Hom}(C(\mu), C(\lambda))\to 
{\rm Hom}(\La, C(\lambda)) \to {\rm Hom}(C(\mu q^{-1}, C(\lambda))$$
$$\to {\rm Ext}^1(C(\mu), C(\lambda))\to 0
$$
With the assumptions, the first and the third term are 1-dimensional. As well ${\rm Hom}(\La, C(\lambda) $ is 2-dimensional, and
hence the ext space is zero.

\bigskip

\begin{cor} Let $M=C(\lambda)$, then ${\rm Ext}^k(M, M)=0$ for $k\geq 2$. Hence $M$ is ext finite and not projective.
\end{cor}

{\it Proof } \ We have ${\rm Ext}^k(M, M) \cong {\rm Ext}^1(\Omega^{k-1}(M), M) = {\rm Ext}^1(C(q^{-k+1}\lambda), C(\lambda)) =0$.

\bigskip

\subsection{Double Nakayama algebras}

We consider algebras of the form 
$A = A(t) = K\cQ/I$ where 
$K\cQ$ is the path algebra of a quiver 
of the form
$$
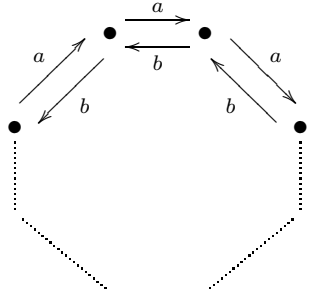
\xymatrix{
& \bullet\ar@<1ex>[r]^a\ar@<1ex>[dl]^{b} & \bullet\ar@<1ex>[l]^{b}\ar@<1ex>[dr]^a & \\
\bullet\ar@<1ex>[ur]^a\ar@{.}[d] & & & \bullet\ar@<1ex>[ul]^{b}\ar@{.}[d]\\
\ar@{.}[rd] & & & \ar@{.}[dl] \\
&&&
}$$
We label the vertices by $\bZ_r$ and the arrows are $a_i: i\mapsto i+1$ and
$b_i: i+1\mapsto i$. The ideal $I$ is generated by 
$a_{i+1}a_{i}, \ \ b_ib_{i+1} $
and 
$$b_ia_i + a_{i-1}b_{i-1}  (i\neq 0), \ \ b_0a_0 + ta_{r-1}b_{r-1} $$
where $0\neq t\in K$. 
We call this algebra, and any Morita equivalent version, a Double Nakayama
algebra.

We want to show that if $t$ is not a root of unity then $A$ has non-projective ext finite modules.

\bigskip

Note that for an arrow $a_i: i\to i+1$ we have in the algebra that 
$a_i = e_{i+1}a_ie_i$ where $e_i$ is the idempotent corresponding to 
vertex $i$.

\bigskip

\begin{lem} The algebra $A$ has a subalgebra $\La$ isomorphic to $\La(q)$ where $q^r-t^{-1}=0$, and 
$A$ is projective as a left and right $\La$-module.
\end{lem}

{\it Proof } \ Let $x: := q^ra_0 + q^{r-1}a_1 + q^{r-2}a_2 + \ldots + qa_{r-1}$ and 
$y:= b_0 + b_1 + b_2 + \ldots + b_{r-1}$. One checks that
$xy+qyx=0$ but $xy\neq 0$; and clearly $x^2=0$ and $y^2=0$. 
Take $B$ to be the  
subalgebra with generators $x, y$.  

Consider $A$ as a left $\La$-module, one checks that $A= \oplus_{i=0}^{r-1} \La e_i$ and that 
$A= \oplus_{i=0}^{r-1} e_i\La$. 

\bigskip

\begin{rem}\normalfont 
(1) By the previous Lemma,  the functor $A\otimes_{\La}(-)$ is exact and takes projective modules to projective modules.
In the following we write $A\otimes(-)$ for $A\otimes_{\La}(-)$. 
As well, for any $\La$-module $N$ the module $A\otimes N$ has dimension $r\cdot \dim N$.

(2) We have $xe_i = qa_i = e_{i+1}x$ and $ye_i = b_{i-1} = e_{i-1}y$. Hence
$$e_i(yx) = (yx)e_i = qb_ia_i$$
(3) If the $\La$-module $N$ has no non-zero projective summands then 
$A\otimes N$ has no non-zero projective summands:
More generally, a module of a selfinjective algebra has no non-zero projective summands  
if and only if it is annihilated by the socle of the algebra.

Here, the socle of $A$ is spanned by the elements $b_ia_i$, and for $w\in N$ we have
$$q(b_ia_i)\otimes w = e_i(yx)\otimes w  = e_i\otimes yxw =0$$
since $yx$ is in the socle of $\La$. 

Now let $M= C(\lambda)$, the $\La$-module as in Subsection 4.1. 

(4) By (1) and (3) we have that  
$$\Omega(A\otimes M) \cong A\otimes \Omega(M) = A\otimes C(q^{-1}\lambda).$$
\end{rem}

\begin{lem}
If $r>0$ then  ${\rm Hom}_A(A\otimes C(q^{-r}\lambda), A\otimes M)$ has dimension $r$.
\end{lem}

\medskip

{\it Proof }  By adjointness 
$${\rm Hom}_A(A\otimes C(q^{-r}\lambda), A\otimes M) \cong {\rm Hom}_{\La}(C(q^{-r}\lambda), A\otimes M)
\leqno{(*)}$$
where 
$A\otimes M$ is restricted to $\La$, and we work with the $\Lambda$-homomorphisms.
One checks that 
the $\La$-socle of $A\otimes M$ is equal to $A\otimes {\rm soc} M = {\rm rad}_{\Lambda}(A\otimes M)$
and hence this has dimension $r$. 

The space (*) contains all maps with image in the $\Lambda$-socle of $A\otimes M$ and this has dimension $r$. So we must
show that for $r\neq 0$ there are no other homomorphisms, that is, we have no monomorphism from $C(q^{-r}\lambda)$ to $A\otimes M$ for $r\neq 0$.

\medskip
  Assume there is a monomorphism, then the image is  a  cyclic $\La$-submodule of $A\otimes M$ of dimension two.
So let $\xi$ be  a generator for a cyclic two-dimensional submodule of $A\otimes M$. We may assume that $\xi$ is of the form
$$\xi = \sum_{i\in \bZ_r} c_i(e_i\otimes \zeta)
$$
[if $w\in {\rm soc}(A)$ then $w\otimes \xi=0$. Furthermore if $w\in {\rm rad}(A)$ and $w\otimes \xi$ is in the socle of $A\otimes M$, then it
may be omitted from a cyclic generator.]

We require that $x\xi$ and $y\xi$ are linearly dependent.
By the identities in Remark 4.5,
$$x\xi = \sum_{j\in \bZ_r} c_{j-1} (e_j\otimes x\zeta), \ \ y\xi= \sum_{j\in \bZ_r} c_{j+1}(e_j\otimes y\zeta)
= \sum_{j\in \bZ_r} \lambda c_{j+1}(e_j\otimes x\zeta).
$$
Assume $y\xi = \mu x\xi$ for some scalar $\mu\neq 0$. The set $\{ e_i\otimes x\zeta\}$  is linearly independent, so we must have
$$c_{j-1} \mu = \lambda c_{j-1}  \ (j\in \bZ_r)
$$
So we get
if $r$ is even, $c_j = (\mu^{-1}\lambda)^{r/2}c_j$ for all $j$ and if $r$ is odd, $c_j = (\mu^{-1}\lambda)^{r-1} c_j$ for all $j$.

Hence if there is such element $\xi$ then $\mu = \lambda\cdot \omega$ for 
some root of unity $\omega$.

If $\mu = q^{-r}\lambda$ then $\mu = \lambda\cdot \omega$ only if $r=0$, and our claim is proved. 
$\Box$

\bigskip

\begin{prop} Let $M=C(\lambda)$. Then $A\otimes M$ is ext-finite and not projective.
\end{prop}

\bigskip

{\it Proof } We have seen in the remark that $A\otimes M$ is not projective,  and we have seen $\Omega^r(A\otimes M) \cong A\otimes C(q^{-r}\lambda)$. 
Take the exact sequence
$$0\to A\otimes C(q^{-r-1}\lambda)\to A\otimes \La \to A\otimes C(q^{-r}\lambda)\to 0
$$
and apply the functor ${\rm Hom}_A(-, A\otimes M)$. Then by using adjointness we get the four-term exact sequence
$$0 \to {\rm Hom}_{\La}(C(q^{-r}\lambda), A\otimes M)
\to {\rm Hom}_{\La}(\La, A\otimes M) \to {\rm Hom}_{\La}(C(q^{-r-1}\lambda), A\otimes M
$$
$$\to {\rm Ext}_{\La}^1(C(q^{-r}), A\otimes M) \cong {\rm Ext}^r_A(\Omega^r(A\otimes M), A\otimes M)\to 0
$$
If $r>0$ then the first and the third term of the sequence has dimension $r$. As well, the second term has dimension $2r$ and hence
the fourth term is zero. 
Hence for all $r>0$ we have ${\rm Ext}^1_A(\Omega^r(A\otimes M), A\otimes M) = 0$. This is isomorphic to
${\rm Ext}^{r+1}_A(A\otimes M, A\otimes M)$ and hence $A\otimes M$ is ext-finite.

$\Box$

\medskip

{\sc Karin Erdmann,  Mathematical Institute, University of Oxford,
ROQ, Oxford OX2 6GG, UK}

\medskip

\textit{email:} {\tt  erdmann@maths.ox.ac.uk}

\end{document}